\documentclass[12pt,leqno]{amsart}
\usepackage{amsfonts,amsmath,amssymb,amsthm, geometry}
\usepackage{graphicx} 
\usepackage{srcltx}
\usepackage[all]{xy}

\usepackage[usenames,dvipsnames]{pstricks}
\usepackage{epsfig}
\usepackage{pst-grad} 
\usepackage{pst-plot} 

\newtheorem{theorem}{Theorem}[section]
\newtheorem{corollary}[theorem]{Corollary}

\newtheorem{proposition}[theorem]{Proposition}
\theoremstyle{definition}
\newtheorem{definition}[theorem]{Definition}
\theoremstyle{remark}
\newtheorem{remark}[theorem]{\sc Remark}
\newtheorem{example}[theorem]{\sc Example}

\newtheorem{note}[theorem]{\sc Note}

\newcommand{\diffeo}{{\rm{diffeo}}}
\newcommand{\codim}{{\rm{codim\hspace{2pt}}}}

\newcommand{\proj}{{\rm{proj}}}

\newcommand{\Sing}{{\rm{Sing\hspace{2pt}}}}

\newcommand{\rank}{{\rm{rank\hspace{2pt}}}}

\newcommand{\re}{{\rm{Re\hspace{2pt}}}}

\newcommand{\im}{{\rm{Im\hspace{2pt}}}}

\newcommand{\cl}{{\rm{closure}}}
\newcommand{\grad}{\mathop{\rm{grad}}}
\newcommand{\ity}{{\infty}}

\newcommand{\e}{\varepsilon}
\newcommand{\m}{\setminus}
\newcommand{\fin}{\hspace*{\fill}$\square$\vspace*{2mm}}

\newcommand{\cS}{{\mathcal S}}

\newcommand{\cW}{{\mathcal W}}

\newcommand{\cN}{{\mathcal N}}

\newcommand{\bR}{{\mathbb R}}
\newcommand{\bC}{{\mathbb C}}

\newcommand{\bN}{{\mathbb N}}

\begin{document}

\title[Open book structures from real mappings]{Singular open book structures from real mappings}

%

\author{\sc R. N. Ara\'ujo dos Santos}
\address{ICMC,
Universidade de S\~ao Paulo,  Av. Trabalhador S\~ao-Carlense, 400 -
CP Box 668, 13560-970 S\~ao Carlos, S\~ao Paulo,  Brazil}
\email{rnonato@icmc.usp.br}

\author{Y. Chen}
\address{Math\'ematiques, UMR-CNRS 8524, Universit\'e Lille 1,
59655 Villeneuve d'Ascq, France.}
\email{Ying.Chen@math.univ-lille1.fr}

\author{M. Tib\u ar}
\address{Math\'ematiques, UMR-CNRS 8524, Universit\'e Lille 1,
59655 Villeneuve d'Ascq, France.}
\email{tibar@math.univ-lille1.fr}



\thanks{R.N. Ara\'ujo dos Santos acknowledges support from  the Brazilian CNPq, PQ-2, proc. num. 305183/2009-5, and FAPESP proc. num. 2009/14383-3.
M. Tib\u ar acknowledges support from the French \textit{Agence Nationale de la Recherche}, grant ANR-08-JCJC-0118-01. 
}


\begin{abstract}
We define open book structures with singular bindings. Starting from an extension of Milnor's results on local fibrations for germs with nonisolated singularity, we find
classes of genuine real analytic mappings which yield such open book structures. 
\end{abstract}

\maketitle


\section{Introduction and main results}\label{intro}



Let $\psi : (\bR^m, 0) \to (\bR^p, 0)$ be some analytic mapping germ, $m > p \ge 2$, and let $V := \psi^{-1}(0)$. The isolated singularity case, i.e. $\Sing \psi = \{0\}$, has been considered by Milnor; his result \cite[Theorem 11.2]{Mi}  tells that there exists a mapping
\begin{equation}\label{eq:milnormapping1}
 S^{m-1}_\e \m K_\e \to S_\delta^{p-1}
\end{equation}
which is a locally trivial fibration and its diffeomorphism type is independent of the radius $\e > 0$ small enough and of $0 <\delta \ll \e$, where $K_\e := V\cap S^{m-1}_\e$. It turns out that the sphere $S^{m-1}_\e$ is moreover endowed with a {\it higher open book structure} with binding $K_\e$, as defined in \cite[Definition 2.1]{dST1}. This is due to the trivial fibration structure induced by $\psi$ in some small neighbourhood $N$ of $K_\e$.   The condition $V  \cap \Sing \psi \subset \{0\}$ used in \cite{dST1} is the most general one under which higher open book structures with smooth  binding $K_\e$ exist. 
Classical open books, i.e. for which the binding is smooth and has codimension 2, are frequent in the literature, they are known under various other names like \textit{fibered links} or \textit{Neuwirth-Stallings pairs} \cite{Lo} or even \textit{Lefschetz pencils}, see e.g. \cite{Wi}.

In this paper we investigate the situation when $\psi^{-1}(0)$ contains nonisolated singularities, thus  the link $K_\e$ is no more a manifold. Several recent papers considered this situation with respect to the existence of the Milnor fibration \eqref{eq:milnormapping1}, see e.g. \cite{PS}, \cite{dST0}, \cite{CSS}.  

We want to give here a precise meaning of ``singular open book structures'', therefore we introduce the following definition:

\begin{definition}\label{d:booksing}
 We say that the pair $(K, \theta)$ is a \textit{higher open book structure with singular binding}  on an analytic manifold $M$ of dimension $m-1 \ge p \ge 2$ if
 $K\subset M$ is a singular real subvariety of codimension $p$ and $\theta : M\setminus K \to S^{p-1}_1$ is a locally trivial smooth fibration such that $K$ admits a neighbourhood $N$ for which the restriction
$\theta_{|N\setminus K}$ is the composition $N\setminus K  \stackrel{h}{\to} B^p \setminus \{ 0\} \stackrel{s/\|s\|}{\to} S^{p-1}_1$ where $h$ is a locally trivial fibration. 

 One says that the \textit{singular fibered link} $K$ is the \textit{binding}  and that the (closures of) the fibers of $\theta$ are the \textit{pages} of the \textit{open book}.
 \end{definition}
From the above definition it also follows that $\theta$ is surjective. 

  In the classical case of open books one has $p=2$, $K\subset M$
is a 2-codimensional submanifold which admits a neighbourhood $N$ diffeomorphic to $B^2 \times K$ 
for which $K$ is identified with $\{ 0\} \times K$ and the restriction  
$\theta_{|N\setminus K}$ is the following composition with the natural projections $N\setminus K  \stackrel{\diffeo}{\simeq} (B^2 \setminus \{ 0\}) \times K \stackrel{\proj}{\to} B^2 \setminus \{ 0\} \to S^1$.  

Therefore our new definition of ``open book with singular binding'' preserves the two main aspects of the classical definitions, namely:

\begin{enumerate}
\item[\bf{(a).}]  the complement of the link fibers over the unit sphere, with codimension $p-1$ fibres.
 \item[\bf{(b).}]  this fibration has a regularity property in the neighbourhood of the link, namely it is induced by a more refined fibration with codimension $p$ fibres.
\end{enumerate}
These two properties turn out to be equally important and independent whenever $V$ contains nonisolated singularities. Our proofs will therefore contain two parts, corresponding to demonstrating property (a) and  property (b). 

Whenever $\psi$ is the pair $(\re f, \im f)$ associated to a holomorphic function germ $f : (\bC^n,0) \to (\bC,0)$,  one has $\Sing \psi \subset V= f^{-1}(0)$ and Milnor \cite[\S 4]{Mi} had proved that the natural mapping:
\begin{equation}\label{eq:holom}
 \frac{f}{|f |} : S_\e^{2n-1}\m K_\e \to S_1^{1}
\end{equation}
is itself a $C^\infty$ locally trivial fibration. This is not enough to provide an open book decomposition whenever $V$ contains nonisolated singularities, since we also need property (b).  Milnor did not show it himself; several years later,  Hamm and L\^e D.T. \cite{HL} proved the nontrivial fact that holomorphic functions have the Thom $(a_f)$-property along a certain stratification of $V$. Then property (b) follows from this.
 
If we consider a holomorphic mapping $\bC^n \to \bC^p$ with $p>1$ then there might not exist fibrations like \eqref{eq:holom}\footnote{in case of an ICIS one does not have the fibration \eqref{eq:holom}.} or the Thom regularity may fail\footnote{e.g. like shown by an example in \cite{HL}.} along $V$.
In the purely real analytic setting, the mapping $\psi$ may have neither property (a) nor property (b). 

\bigskip

We formulate in \S \ref{s:milnor} an extension of \cite[Theorem 11.2]{Mi} which uses our Definition \ref{d:booksing}, the proof of which is derived from carefully revisiting Milnor's construction. This is further used as a basis for our main results. In order to state them, we need the following definitions.

\begin{definition}\label{d:M}
Let $U\subset \bR^m$ be an open set and let $\rho : U \to \bR_{\ge 0}$ be a proper analytic function. 
We say that the set of \textit{$\rho$-nonregular points}\footnote{sometimes called \textit{Milnor set}, e.g. in \cite{NZ}, the name ``$\rho$-regularity'' has been introduced in \cite{Ti-cras}.} of an analytic mapping  $\Psi :U \to \bR^p$ is the set of non-transversality between $\rho$ and $\Psi$, i.e.
 $M(\Psi) := \{ x\in U \mid \rho \not\pitchfork_x \Psi\}$.
 
 Similarly, the set of $\rho$-nonregular points of the mapping $\frac{\Psi}{\|\Psi\|} :U\m V \to S^{p-1}_1$ is the set: $M(\frac{\Psi}{\|\Psi\|}) :=  \cl \{ x\in U \m V \mid \rho \not\pitchfork_x \frac{\Psi}{\|\Psi\|}\}$.
 \end{definition}

Let us remind that, by definition, two mappings on an open $U\subset \bR^m$ are called transversal at some point $x\in U$ iff they are both non-singular at $x$ and their well-defined tangent spaces at $x$ are transversal in $\bR^m$. In particular, the singular locus of each of the mappings is included in their non-transversality locus.

The \textit{$\rho$-regularity} of Definition \ref{d:M} is a basic tool, used by many authors (Thom, Milnor, Mather, Looijenga, Bekka etc) in the local stratified setting as well as at infinity, like in e.g. \cite{NZ}, [Ti1,2,3] in order to produce locally trivial fibrations. See also \cite{CSS}. Thom has called such a $\rho$ ``fonction tapissante''. As also done in \cite{dST0}, in this paper we shall use as function $\rho$ the Euclidean distance function and its open balls and spheres of radius $\e$, denoted by $B_\e$ and $S_\e$ respectively. The open set $U$ of Definition  \ref{d:M} will be small ball $B_\e$, unless otherwise stated.
 
\begin{theorem}\label{t:isolcrt}
 Let $\psi : (\bR^m, 0) \to (\bR^p, 0)$ be an analytic mapping germ, $m> p\ge 2$,  such that $\codim V = p$.
  Suppose that:
\begin{equation}\label{eq:main}
\overline{M(\psi)\m V} \cap  V = \{ 0\}. 
\end{equation}

If $M(\frac{\psi}{\|\psi\|}) = \emptyset$ then $(K_\e, \frac{\psi}{\|\psi \|})$ is an open book structure with singular binding on $S^{m-1}_\e$,   independent (up to isotopies) of $\e >0$ small enough.
\end{theorem}

 Condition \eqref{eq:main} allows nonisolated singularities, more precisely it implies $M(\psi) = A \cup B$ where $A\subset V$ and $B \cap V \subset \{ 0\}$, and both $A$ and $B$ may be of positive dimension. We also have $\Sing \psi \subset M(\psi)$. Therefore Theorem \ref{t:isolcrt} represents a simultaneous extension of \cite[Theorem 2.2]{dST1}, where the singular locus was of type B, and of \cite[Proposition 5.3]{dST0} or \cite[Theorem 5.3]{CSS}\footnote{compare the result \cite[Proposition 5.3 and Remark 5.4]{dST0} to \cite[Theorem 5.3 points (1),(2) and Remark 5.7]{CSS}. Then compare the comments \cite[p. 423, lines 5-8]{CSS} to the reality of \cite[Section 5]{dST0}.}, where the singular locus was of type A.
We discuss in \S \ref{s:strong} other particular cases and the relations to \cite{PS}.
 
It turns out that the condition\footnote{Massey \cite{Ma} used condition \eqref{eq:main} in conjunction with $\Sing \psi \subset V$ in order to get a full tube fibration \eqref{eq:tube}.} \eqref{eq:main}  insures the existence of the locally trivial fibration $N\setminus K  \stackrel{h}{\to} B^p \setminus \{ 0\}$ from Definition \ref{d:booksing} and is implied by the Thom regularity at $V$. We shall discuss this relation and other criteria in \S \ref{s:classicandnew}. We show by Example \ref{e:thom} how to check this condition directly.

\smallskip
In order to produce a new class of higher dimensional purely real examples, we use the theory of \textit{mixed functions}, i.e. real analytic mappings $\bR^{2n} \simeq \bC^n \to \bC \simeq \bR^2$, recently developed by Mutsuo Oka (see \cite{oka1, oka2} and our footnote at \S \ref{ss:polar}). The necessary definitions are given in \S \ref{ss:polar}. We only point out here that for mixed functions, unlike the holomorphic ones,  the notion of ``polar weighted-homogeneous'' is different and independent from ``radial weighted-homogeneous''. The later is discussed in \S \ref{s:strong} and is one of the hypotheses of Corollary \ref{c:radial}. The former occurs in the following statement (discussed and proved in \S \ref{ss:polar}). 

 \begin{theorem}\label{t:homogen}
\sloppy Let $f:\mathbb{C}^{n}\rightarrow\mathbb{C}$ be a non-constant mixed polynomial which is polar weighted-homogeneous, $n\ge 2$, such that $\codim_\bR V = 2$.
 Then $(K_\e, \frac{f}{\|f \|})$ is an open book structure with singular binding on $S^{2n-1}_\e$,
   independent (up to isotopies) of $\e >0$ small enough. 
\end{theorem}

We finally present some other new classes of examples in \S \ref{s:classicandnew}, one of them by using a Thom-Sebastiani type statement, Proposition \ref{p:exempl}, which represents a new result in the real context.  

\begin{note}\label{n:singularspaceX}
 One may work in a slightly more general setting than Definition \ref{d:booksing}, as follows. Instead of a manifold $M$, let $M$ be a connected compact real analytic set with $\Sing M \subset K$. In this paper $M$ is the link of $(\bR^m, 0)$, hence a sphere, but we may also consider a real analytic germ $(X, 0) \subset (\bR^m, 0)$ with connected link $M$ with respect to some distance function which is not necessarily the Euclidean one, and mappings $\psi : (X, 0)\to (\bR^p, 0)$ such that $\Sing X \subset \psi^{-1}(0)$.  Then one has to modify Milnor's proof  \S \ref{ss:blow} such that the vector field along which one blows the tube to the sphere is tangent to $X$. See also Remark \ref{r:pichonseade}. 
 
 Even more generally, if one considers any singular stratified space $X$ and its link $M$, then Milnor's method reviewed in \S \ref{ss:milnor} still works, i.e. extends (by using classical technical devices of ``radial vector fields'') for stratified transversality and stratified vector fields. This yields what one could call \textit{open book structures with singular pages}.
\end{note}

\section{Milnor's method and open book structures}\label{s:milnor}

In the real analytic setting, Milnor observed that the fibration (\ref{eq:milnormapping1}) may  not be induced by the mapping $\theta = \frac{\psi}{\|\psi \|}$. He gave an example \cite[p. 99]{Mi} of $\psi$ with isolated singularity,  $m=p=2$ and $K_\e= \emptyset$, showing that the mapping 
 \begin{equation}\label{eq:milnormapping2}
 \frac{\psi}{\|\psi \|} : S_\e^{m-1}\m K_\e \to S_1^{1}
\end{equation}
 is not a submersion, hence not a locally trivial fibration.

 In case of isolated singularity and $p=2$, several authors obtained sufficient conditions under which (\ref{eq:milnormapping2}) is a fibration \cite{Ja1, Ja2, RS, RSV} and provided examples showing that the class of real mapping germs $\psi$ with isolated singularity satisfying them enlarges the class of holomorphic functions $f$. 
 For a more general setting $\Sing \psi \cap V \subset \{ 0\}$ and any $m\ge p \ge 2$ an existence criterion has been given in \cite{dST1} and shown to be more general than the previous ones from the literature.

The \textit{$\rho$-regularity} of Definition \ref{d:M}
expresses the transversality of the fibres of a mapping $\psi$ to the levels of $\rho$. It is a basic tool, used by many authors (Milnor, Mather, Looijenga, Bekka etc) in the local stratified setting as well as at infinity, like in e.g. \cite[Ti1,2,3]{NZ}, in order to produce locally trivial fibrations.  As also done in \cite{dST0}, in this paper we shall use as function $\rho$ the Euclidean distance function and its open balls and spheres of radius $\e$, currently denoted by $B_\e$ and $S_\e$ respectively. 

From Definition \ref{d:M} it follows that $M(\psi)$ is a relatively closed analytic set containing the singular set $\Sing \psi$. As for $M(\frac{\psi}{\|\psi\|})$, it is by definition closed but does not necessarily include $\Sing \psi$. We nevertheless have $M(\frac{\psi}{\|\psi\|})\m V \subset M(\psi)\m V$, since $\rho \pitchfork_x \psi$ implies $\rho \pitchfork_x \frac{\psi}{\|\psi\|}$ for $x\not\in V$. In the following we tacitly conceive these non-regularity sets as \textit{set germs at the origin}.

\begin{theorem}\label{t:milnor} \rm (after \cite{Mi}) \it \\
 Let $\psi : (\bR^m, 0) \to (\bR^p, 0)$, $m>p \ge 2$,  be an analytic mapping germ with $\Sing \psi \subset V$, $\codim_\bR V = p$ and satisfying the condition \eqref{eq:main}. Then there exists a higher open book structure with singular binding $(K_\e, \theta)$ on $S^{m-1}_\e$, which
  is independent of $\e >0$ small enough, up to isotopies.
\end{theorem}

\subsection{Revisiting  Milnor's method}\label{ss:milnor}
In case of mapping germs $\psi :(\bR^m , 0) \to (\bR^p, 0)$, in
Milnor's proof \cite[p. 97-99, Theorem 11.2]{Mi} of the existence of a locally trivial fibration
(\ref{eq:milnormapping1})
 one may distinguish two key parts, already pointed out in \cite[\S 5]{dST0}. We shortly review them below in order to show how they apply to the more general situation displayed in Theorem \ref{t:milnor}.

\subsubsection{Existence of the tube fibration}
Assume that there exists $\e_0 >0$ such that the mapping:
\begin{equation}\label{eq:tube}
 \psi_| : \bar B^{m}_{\e} \cap \psi^{-1}(\bar B^{p}_\eta \m \{ 0\}) \to \bar B^{p}_\eta \m \{ 0\}
\end{equation}
 is a surjective locally trivial $C^\ity$-fibration, for all $0< \e \le \e_0$ and $0< \eta= \eta(\e) \ll \e$.
 We shall call it \textit{the tube fibration}. 

 In the case $\Sing \psi = \{0\}$ considered by Milnor, $V$ is transverse to all small enough spheres and therefore any such sphere is also transverse to all nearby fibres, which are moreover non-singular. By Ehresmann's theorem for manifolds with boundary, one concludes the existence of the tube fibration.
 
\subsubsection{Existence of a fibration in the exterior of the tube} \label{ss:blow}
This is Milnor's proof of ``inflating'' the empty tube $\bar B^{m}_{\e} \cap \psi^{-1}(S^{p-1}_\eta)$ to the sphere $S^{m-1}_{\e}\m K_\e$.  Milnor explains in \cite[p. 99]{Mi} that, given a real analytic mapping $\psi$, one may construct, like in \cite[Lemmas 11.3, 5.9, 5.10]{Mi}, a nowhere zero $C^\ity$-vector field $v(x)$ on $\bar B^{m}_{\e}\m \psi^{-1}(B^{p}_\eta)$ satisfying the following two conditions:
$\left\langle  x , v(x)\right\rangle > 0$ and
$\left\langle \grad \|\psi (x)\|^2, v(x) \right \rangle >0$.
The first condition says that $v(x)$ is transverse to all small spheres and points outwards,  and the second condition says that the mapping $\|\psi (x)\|^2$ increases along the flow.
This vector field may be integrated and produces a diffeomorphism:
 \begin{equation}\label{eq:diff}
 \gamma : \bar B^{m}_{\e} \cap \psi^{-1}(S^{p-1}_\eta)\to  S^{m-1}_\e \m \psi^{-1}(B^{p}_\eta).
\end{equation}

This procedure combines the position vector field  $x$ with the gradient vector field $\grad \| \psi \|^2$  and works if the latter is nowhere zero in the neighbourhood of $\psi^{-1}( B^{p}_\eta)$, if the two vector fields never point in the opposite directions \cite[Cor. 3.4]{Mi}, and if the empty tube $\bar B^{m}_{\e} \cap \psi^{-1}(S^{p-1}_\eta)$ is a manifold with boundary.  Milnor uses this construction for holomorphic functions $f$, where $\Sing f \subset V$ and for real mappings with $\Sing \psi = \{ 0\}$. Milnor's technique holds locally at $0$ in the more 
general setting $\Sing \psi \subset V$.

 \subsubsection{Conclusions}
If the tube fibration (\ref{eq:tube}) exists, then its restriction to the empty tube:
\begin{equation}\label{eq:emptytube}
 \psi_| : \bar B^{m}_{\e} \cap \psi^{-1}(S^{p-1}_\eta) \to  S^{p-1}_\eta
\end{equation}
 is a locally trivial $C^\ity$-fibration too.
 
If the inflating procedure works, then the diffeomorphism (\ref{eq:diff}) induces a mapping
$\mu : S^{m-1}_\e \m \psi^{-1}( B^{p}_\eta) \to  S^{p-1}_\eta$
 which is a locally trivial fibration and coincides with $\psi$ on $S^{m-1}_\e \cap \psi^{-1}( S^{p-1}_\eta)$. This may be composed with the mapping $\frac{s}{\| s\|} : S^{p-1}_\eta \to S^{p-1}_1$ and yields a locally trivial fibration
\begin{equation}\label{eq:sphere}
 \mu' : S^{m-1}_\e \m \psi^{-1}( B^{p}_\eta) \to  S^{p-1}_1.
\end{equation}

\subsection{Proof of Theorem \ref{t:milnor}}\label{ss:endproof}
 The hypothesis $\overline{M(\psi)\m V} \cap  V = \{ 0\}$ is equivalent to the existence of a neighbourhood $\cN$ of $V \m \{ 0\}$ such that $M(\psi) \cap \cN \m V = \emptyset$ (see Figure). Together with $\Sing \psi \subset V$, these imply that for any $\e > 0$ small enough, there exists some positive $\eta \ll \e$ such that the mapping (\ref{eq:tube})
is a proper submersion over the pointed open ball $B^p\m \{ 0\}$. Indeed, the properness follows since the restriction of $\psi$ to $\bar B^m_{\epsilon} \m V$ is proper, and ``submersion'' is a consequence of the condition $\Sing \psi \subset V$ and of the transversality of the fibres to the boundary $S^{m-1}_\e = \partial \bar B^m_{\epsilon}$. It then follows that this restriction is also surjective. Just like Milnor did,  we may now apply  Ehresmann's theorem \cite{Eh, Wo} for manifolds with boundary to conclude to the existence of the locally trivial fibration (\ref{eq:tube}).\footnote{this observation is also contained in \cite[Theorem 4.4]{Ma}.}

   To show now that the fibration (\ref{eq:sphere}) extends to an open book structure on $S^{m-1}_\e$ one must produce the map $\theta$, cf Definition \ref{d:booksing}. The fibration (\ref{eq:sphere}) may be glued along $S^{m-1}_\e \cap \psi^{-1}( S^{p-1}_\eta)$ to the locally trivial fibration $S^{m-1}_\e \cap \psi^{-1}(\bar B^{p}_\eta \m \{ 0\}) \to \bar B^{p}_\eta \m \{ 0\}$
 composed with the mapping $\frac{s}{\| s\|} : \bar B^{p}_\eta \m \{ 0\} \to S^{p-1}_1$ 
since their restrictions to this boundary coincide. This glueing may be done in the $C^\ity$ category and produces 
a locally trivial $C^ \ity$-fibration. We then define $\theta$ to be the result of the glueing of $\mu'$ with  $\frac{s}{\| s\|}\circ \psi_|  :  S^{m-1}_{\e} \cap \psi^{-1}(\bar B^{p}_\eta \m \{ 0\}) \to \bar B^{p}_\eta \m \{ 0\} \to S^{p-1}_1$, and get that $(K_\e, \theta)$ is an open book decomposition of $S^{m-1}_\e$.
\fin

\vspace{0.5cm}
\begin{center}
\scalebox{1} 
{
\begin{pspicture}(0,-4.3940544)(8.961875,4.3792186)
\definecolor{color282}{rgb}{0.0,0.2,1.0}
\definecolor{color420}{rgb}{0.0,0.2,0.8}
\pscircle[linewidth=0.04,dimen=outer](3.68,-0.49921876){2.7}
\psdots[dotsize=0.12](3.62,-0.49921876)
\pscustom[linewidth=0.04,linecolor=blue]
{
\newpath
\moveto(3.6,-0.47921875)
\lineto(3.7,-0.44921875)
\curveto(3.75,-0.43421876)(3.865,-0.38421875)(3.93,-0.34921876)
\curveto(3.995,-0.31421876)(4.13,-0.16921875)(4.2,-0.05921875)
\curveto(4.27,0.05078125)(4.365,0.22578125)(4.39,0.29078126)
\curveto(4.415,0.35578126)(4.5,0.5107812)(4.56,0.60078126)
\curveto(4.62,0.69078124)(4.72,0.8957813)(4.76,1.0107813)
\curveto(4.8,1.1257813)(4.86,1.4857812)(4.88,1.7307812)
\curveto(4.9,1.9757812)(4.925,2.3557813)(4.93,2.4907813)
\curveto(4.935,2.6257813)(4.95,2.8257813)(4.96,2.8907812)
\curveto(4.97,2.9557812)(4.98,3.0257812)(4.98,3.0407813)
}
\pscustom[linewidth=0.04,linecolor=blue]
{
\newpath
\moveto(3.6,-0.45921874)
\lineto(3.68,-0.57921875)
\curveto(3.72,-0.63921875)(3.795,-0.76921874)(3.83,-0.83921874)
\curveto(3.865,-0.9092187)(3.915,-1.0192188)(3.93,-1.0592188)
\curveto(3.945,-1.0992187)(3.995,-1.2292187)(4.03,-1.3192188)
\curveto(4.065,-1.4092188)(4.135,-1.6192187)(4.17,-1.7392187)
\curveto(4.205,-1.8592187)(4.255,-2.0492187)(4.27,-2.1192188)
\curveto(4.285,-2.1892188)(4.305,-2.3442187)(4.31,-2.4292188)
\curveto(4.315,-2.5142188)(4.32,-2.6592188)(4.32,-2.7192187)
\curveto(4.32,-2.7792187)(4.32,-2.8842187)(4.32,-2.9292188)
\curveto(4.32,-2.9742188)(4.31,-3.1542187)(4.3,-3.2892187)
\curveto(4.29,-3.4242187)(4.275,-3.5942187)(4.27,-3.6292188)
\curveto(4.265,-3.6642187)(4.26,-3.7692187)(4.26,-3.8392189)
\curveto(4.26,-3.9092188)(4.265,-4.0942187)(4.27,-4.209219)
}
\pscustom[linewidth=0.04,shadow=true,shadowangle=-45.0]
{
\newpath
\moveto(3.58,-0.43921876)
\lineto(3.64,-0.41921875)
\curveto(3.67,-0.40921876)(3.715,-0.38921875)(3.73,-0.37921876)
\curveto(3.745,-0.36921874)(3.78,-0.35421875)(3.8,-0.34921876)
\curveto(3.82,-0.34421876)(3.86,-0.32921875)(3.88,-0.31921875)
\curveto(3.9,-0.30921876)(3.935,-0.28421876)(3.95,-0.26921874)
\curveto(3.965,-0.25421876)(3.995,-0.21921875)(4.01,-0.19921875)
\curveto(4.025,-0.17921875)(4.05,-0.13421875)(4.06,-0.10921875)
\curveto(4.07,-0.08421875)(4.1,-0.05421875)(4.12,-0.04921875)
\curveto(4.14,-0.04421875)(4.175,-0.01921875)(4.19,0.0)
\curveto(4.205,0.02078125)(4.225,0.06078125)(4.23,0.08078125)
\curveto(4.235,0.10078125)(4.25,0.14578125)(4.26,0.17078125)
\curveto(4.27,0.19578125)(4.29,0.24078125)(4.3,0.26078126)
\curveto(4.31,0.28078124)(4.34,0.31578124)(4.36,0.33078125)
\curveto(4.38,0.34578124)(4.41,0.37578124)(4.42,0.39078125)
\curveto(4.43,0.40578124)(4.45,0.45078126)(4.46,0.48078126)
\curveto(4.47,0.5107812)(4.49,0.5757812)(4.5,0.61078125)
\curveto(4.51,0.6457813)(4.525,0.7157813)(4.53,0.75078124)
\curveto(4.535,0.78578126)(4.555,0.85578126)(4.57,0.8907812)
\curveto(4.585,0.92578125)(4.62,0.97078127)(4.64,0.98078126)
\curveto(4.66,0.99078125)(4.69,1.0207813)(4.7,1.0407813)
\curveto(4.71,1.0607812)(4.725,1.1107812)(4.73,1.1407813)
\curveto(4.735,1.1707813)(4.76,1.2107812)(4.78,1.2207812)
\curveto(4.8,1.2307812)(4.83,1.2607813)(4.84,1.2807813)
\curveto(4.85,1.3007812)(4.86,1.3707813)(4.86,1.4207813)
\curveto(4.86,1.4707812)(4.86,1.5557812)(4.86,1.5907812)
\curveto(4.86,1.6257813)(4.86,1.7257812)(4.86,1.7907813)
\curveto(4.86,1.8557812)(4.855,1.9507812)(4.85,1.9807812)
\curveto(4.845,2.0107813)(4.84,2.0857813)(4.84,2.1307812)
\curveto(4.84,2.1757812)(4.845,2.2507813)(4.85,2.2807813)
\curveto(4.855,2.3107812)(4.865,2.3657813)(4.87,2.3907812)
\curveto(4.875,2.4157813)(4.88,2.4757812)(4.88,2.5107813)
\curveto(4.88,2.5457811)(4.89,2.6157813)(4.9,2.6507812)
\curveto(4.91,2.6857812)(4.92,2.7457812)(4.92,2.7707813)
\curveto(4.92,2.7957811)(4.925,2.8457813)(4.93,2.8707812)
\curveto(4.935,2.8957813)(4.94,2.9457812)(4.94,2.9707813)
\curveto(4.94,2.9957812)(4.95,3.0357811)(4.96,3.0507812)
\curveto(4.97,3.0657814)(4.98,3.0857813)(4.98,3.1007812)
\openshadow
}
\pscustom[linewidth=0.04,shadow=true,shadowangle=-45.0]
{
\newpath
\moveto(3.6,-0.49921876)
\lineto(3.63,-0.55921876)
\curveto(3.645,-0.58921874)(3.675,-0.63421875)(3.69,-0.64921874)
\curveto(3.705,-0.6642187)(3.735,-0.69421875)(3.75,-0.70921874)
\curveto(3.765,-0.7242187)(3.795,-0.76421875)(3.81,-0.7892187)
\curveto(3.825,-0.81421876)(3.845,-0.90421873)(3.85,-0.96921873)
\curveto(3.855,-1.0342188)(3.865,-1.1342187)(3.87,-1.1692188)
\curveto(3.875,-1.2042187)(3.885,-1.2592187)(3.89,-1.2792188)
\curveto(3.895,-1.2992188)(3.92,-1.3242188)(3.94,-1.3292187)
\curveto(3.96,-1.3342187)(3.995,-1.3642187)(4.01,-1.3892188)
\curveto(4.025,-1.4142188)(4.05,-1.4642187)(4.06,-1.4892187)
\curveto(4.07,-1.5142188)(4.105,-1.6242187)(4.13,-1.7092187)
\curveto(4.155,-1.7942188)(4.19,-1.9042188)(4.2,-1.9292188)
\curveto(4.21,-1.9542187)(4.225,-2.0592186)(4.23,-2.1392188)
\curveto(4.235,-2.2192187)(4.245,-2.3242188)(4.25,-2.3492188)
\curveto(4.255,-2.3742187)(4.27,-2.4342186)(4.28,-2.4692187)
\curveto(4.29,-2.5042188)(4.305,-2.5742188)(4.31,-2.6092188)
\curveto(4.315,-2.6442187)(4.32,-2.7042189)(4.32,-2.7292187)
\curveto(4.32,-2.7542188)(4.32,-2.8242188)(4.32,-2.8692188)
\curveto(4.32,-2.9142187)(4.32,-2.9842188)(4.32,-3.0092187)
\curveto(4.32,-3.0342188)(4.315,-3.0792189)(4.31,-3.0992188)
\curveto(4.305,-3.1192188)(4.295,-3.1642187)(4.29,-3.1892188)
\curveto(4.285,-3.2142189)(4.28,-3.2842188)(4.28,-3.3292189)
\curveto(4.28,-3.3742187)(4.275,-3.4492188)(4.27,-3.4792187)
\curveto(4.265,-3.5092187)(4.25,-3.5992188)(4.24,-3.6592188)
\curveto(4.23,-3.7192187)(4.215,-3.8292189)(4.21,-3.8792188)
\curveto(4.205,-3.9292188)(4.205,-4.0242186)(4.21,-4.0692186)
\curveto(4.215,-4.1142187)(4.225,-4.184219)(4.23,-4.209219)
\curveto(4.235,-4.2342186)(4.24,-4.269219)(4.24,-4.2992187)
\openshadow
}
\pscustom[linewidth=0.04,linecolor=blue]
{
\newpath
\moveto(3.56,-0.47921875)
\lineto(3.48,-0.40921876)
\curveto(3.44,-0.37421876)(3.365,-0.30921876)(3.33,-0.27921876)
\curveto(3.295,-0.24921875)(3.18,-0.17421874)(3.1,-0.12921876)
\curveto(3.02,-0.08421875)(2.85,0.0)(2.76,0.04078125)
\curveto(2.67,0.08078125)(2.5,0.12078125)(2.42,0.12078125)
\curveto(2.34,0.12078125)(2.19,0.10578125)(2.12,0.09078125)
\curveto(2.05,0.07578125)(1.87,0.08078125)(1.76,0.10078125)
\curveto(1.65,0.12078125)(1.465,0.14078125)(1.39,0.14078125)
\curveto(1.315,0.14078125)(1.17,0.16078125)(1.1,0.18078125)
\curveto(1.03,0.20078126)(0.875,0.26078126)(0.79,0.30078125)
\curveto(0.705,0.34078124)(0.555,0.44078124)(0.49,0.50078124)
\curveto(0.425,0.56078124)(0.29,0.61078125)(0.22,0.60078126)
\curveto(0.15,0.5907813)(0.06,0.5807812)(0.0,0.5807812)
}
\pscustom[linewidth=0.04,linecolor=blue]
{
\newpath
\moveto(3.6,-0.47921875)
\lineto(3.64,-0.49921876)
\curveto(3.66,-0.50921875)(3.705,-0.5342187)(3.73,-0.5492188)
\curveto(3.755,-0.56421876)(3.815,-0.59421873)(3.85,-0.6092188)
\curveto(3.885,-0.62421876)(3.97,-0.6692188)(4.02,-0.69921875)
\curveto(4.07,-0.7292187)(4.185,-0.7942188)(4.25,-0.82921875)
\curveto(4.315,-0.8642188)(4.43,-0.9142187)(4.48,-0.92921877)
\curveto(4.53,-0.94421875)(4.675,-0.95421875)(4.77,-0.94921875)
\curveto(4.865,-0.94421875)(5.05,-0.93921876)(5.14,-0.93921876)
\curveto(5.23,-0.93921876)(5.37,-0.93921876)(5.42,-0.93921876)
\curveto(5.47,-0.93921876)(5.6,-0.93921876)(5.68,-0.93921876)
\curveto(5.76,-0.93921876)(5.88,-0.94921875)(5.92,-0.95921874)
\curveto(5.96,-0.96921873)(6.115,-1.0192188)(6.23,-1.0592188)
\curveto(6.345,-1.0992187)(6.525,-1.1792188)(6.59,-1.2192187)
\curveto(6.655,-1.2592187)(6.74,-1.3342187)(6.76,-1.3692187)
\curveto(6.78,-1.4042188)(6.845,-1.4792187)(6.89,-1.5192188)
\curveto(6.935,-1.5592188)(7.005,-1.6192187)(7.03,-1.6392188)
\curveto(7.055,-1.6592188)(7.085,-1.6842188)(7.1,-1.6992188)
}
\pscustom[linewidth=0.04,linestyle=dashed,dash=0.16cm 0.16cm]
{
\newpath
\moveto(3.56,-0.43921876)
\lineto(3.59,-0.40921876)
\curveto(3.605,-0.39421874)(3.63,-0.36421874)(3.64,-0.34921876)
\curveto(3.65,-0.33421874)(3.675,-0.29421875)(3.69,-0.26921874)
\curveto(3.705,-0.24421875)(3.735,-0.19921875)(3.75,-0.17921875)
\curveto(3.765,-0.15921874)(3.795,-0.11921875)(3.81,-0.09921875)
\curveto(3.825,-0.07921875)(3.855,-0.02921875)(3.87,0.0)
\curveto(3.885,0.03078125)(3.91,0.09578125)(3.92,0.13078125)
\curveto(3.93,0.16578124)(3.96,0.23078126)(3.98,0.26078126)
\curveto(4.0,0.29078126)(4.02,0.34578124)(4.02,0.37078124)
\curveto(4.02,0.39578125)(4.045,0.45078126)(4.07,0.48078126)
\curveto(4.095,0.5107812)(4.135,0.63078123)(4.15,0.72078127)
\curveto(4.165,0.81078124)(4.19,0.94078124)(4.2,0.98078126)
\curveto(4.21,1.0207813)(4.24,1.1107812)(4.26,1.1607813)
\curveto(4.28,1.2107812)(4.315,1.3157812)(4.33,1.3707813)
\curveto(4.345,1.4257812)(4.37,1.5607812)(4.38,1.6407813)
\curveto(4.39,1.7207812)(4.405,1.9057813)(4.41,2.0107813)
\curveto(4.415,2.1157813)(4.42,2.3107812)(4.42,2.4007812)
\curveto(4.42,2.4907813)(4.425,2.6507812)(4.43,2.7207813)
\curveto(4.435,2.7907813)(4.44,2.9107811)(4.44,2.9607813)
\curveto(4.44,3.0107813)(4.44,3.0857813)(4.44,3.1107812)
\curveto(4.44,3.1357813)(4.44,3.1657813)(4.44,3.1807814)
}
\pscustom[linewidth=0.04,linestyle=dashed,dash=0.16cm 0.16cm]
{
\newpath
\moveto(3.64,-0.49921876)
\lineto(3.67,-0.47921875)
\curveto(3.685,-0.46921876)(3.735,-0.45921874)(3.77,-0.45921874)
\curveto(3.805,-0.45921874)(3.89,-0.45921874)(3.94,-0.45921874)
\curveto(3.99,-0.45921874)(4.08,-0.45421875)(4.12,-0.44921875)
\curveto(4.16,-0.44421875)(4.22,-0.42421874)(4.24,-0.40921876)
\curveto(4.26,-0.39421874)(4.31,-0.35921875)(4.34,-0.33921874)
\curveto(4.37,-0.31921875)(4.425,-0.28921875)(4.45,-0.27921876)
\curveto(4.475,-0.26921874)(4.52,-0.23421875)(4.54,-0.20921876)
\curveto(4.56,-0.18421875)(4.6,-0.13921875)(4.62,-0.11921875)
\curveto(4.64,-0.09921875)(4.68,-0.03421875)(4.7,0.01078125)
\curveto(4.72,0.05578125)(4.745,0.13078125)(4.75,0.16078125)
\curveto(4.755,0.19078125)(4.765,0.25578126)(4.77,0.29078126)
\curveto(4.775,0.32578126)(4.805,0.37578124)(4.83,0.39078125)
\curveto(4.855,0.40578124)(4.89,0.44078124)(4.9,0.46078125)
\curveto(4.91,0.48078126)(4.935,0.54578125)(4.95,0.5907813)
\curveto(4.965,0.6357812)(4.99,0.74578124)(5.0,0.81078124)
\curveto(5.01,0.87578124)(5.035,0.98578125)(5.05,1.0307813)
\curveto(5.065,1.0757812)(5.09,1.1407813)(5.1,1.1607813)
\curveto(5.11,1.1807812)(5.135,1.2507813)(5.15,1.3007812)
\curveto(5.165,1.3507812)(5.2,1.4257812)(5.22,1.4507812)
\curveto(5.24,1.4757812)(5.26,1.5357813)(5.26,1.5707812)
\curveto(5.26,1.6057812)(5.275,1.6857812)(5.29,1.7307812)
\curveto(5.305,1.7757813)(5.325,1.8707813)(5.33,1.9207813)
\curveto(5.335,1.9707812)(5.36,2.0807812)(5.38,2.1407812)
\curveto(5.4,2.2007813)(5.43,2.3107812)(5.44,2.3607812)
\curveto(5.45,2.4107811)(5.465,2.4857812)(5.47,2.5107813)
\curveto(5.475,2.5357811)(5.485,2.5957813)(5.49,2.6307812)
\curveto(5.495,2.6657813)(5.505,2.7207813)(5.51,2.7407813)
\curveto(5.515,2.7607813)(5.52,2.8007812)(5.52,2.8607812)
}
\pscustom[linewidth=0.04,linestyle=dashed,dash=0.16cm 0.16cm]
{
\newpath
\moveto(3.56,-0.51921874)
\lineto(3.56,-0.63921875)
\curveto(3.56,-0.69921875)(3.56,-0.80421877)(3.56,-0.8492187)
\curveto(3.56,-0.89421874)(3.56,-1.0042187)(3.56,-1.0692188)
\curveto(3.56,-1.1342187)(3.56,-1.2392187)(3.56,-1.2792188)
\curveto(3.56,-1.3192188)(3.575,-1.3892188)(3.59,-1.4192188)
\curveto(3.605,-1.4492188)(3.645,-1.5642188)(3.67,-1.6492188)
\curveto(3.695,-1.7342187)(3.73,-1.8642187)(3.74,-1.9092188)
\curveto(3.75,-1.9542187)(3.775,-2.0542188)(3.79,-2.1092188)
\curveto(3.805,-2.1642187)(3.825,-2.3592188)(3.83,-2.4992187)
\curveto(3.835,-2.6392188)(3.84,-2.8842187)(3.84,-2.9892187)
\curveto(3.84,-3.0942187)(3.845,-3.2542188)(3.85,-3.3092186)
\curveto(3.855,-3.3642187)(3.86,-3.5692186)(3.86,-3.7192187)
\curveto(3.86,-3.8692188)(3.865,-4.084219)(3.87,-4.1492186)
\curveto(3.875,-4.2142186)(3.88,-4.2942185)(3.88,-4.3392186)
}
\pscustom[linewidth=0.04,linestyle=dashed,dash=0.16cm 0.16cm]
{
\newpath
\moveto(3.64,-0.49921876)
\lineto(3.68,-0.5392187)
\curveto(3.7,-0.55921876)(3.735,-0.58921874)(3.75,-0.5992187)
\curveto(3.765,-0.6092188)(3.795,-0.63421875)(3.81,-0.64921874)
\curveto(3.825,-0.6642187)(3.86,-0.69421875)(3.88,-0.70921874)
\curveto(3.9,-0.7242187)(3.935,-0.75921875)(3.95,-0.77921873)
\curveto(3.965,-0.7992188)(4.0,-0.84421873)(4.02,-0.86921877)
\curveto(4.04,-0.89421874)(4.095,-0.9742187)(4.13,-1.0292188)
\curveto(4.165,-1.0842187)(4.23,-1.1892188)(4.26,-1.2392187)
\curveto(4.29,-1.2892188)(4.34,-1.3642187)(4.36,-1.3892188)
\curveto(4.38,-1.4142188)(4.405,-1.5092187)(4.41,-1.5792187)
\curveto(4.415,-1.6492188)(4.455,-1.8342187)(4.49,-1.9492188)
\curveto(4.525,-2.0642188)(4.58,-2.2342188)(4.6,-2.2892187)
\curveto(4.62,-2.3442187)(4.64,-2.4592187)(4.64,-2.5192187)
\curveto(4.64,-2.5792189)(4.65,-2.6842186)(4.66,-2.7292187)
\curveto(4.67,-2.7742188)(4.69,-2.8492188)(4.7,-2.8792188)
\curveto(4.71,-2.9092188)(4.74,-2.9742188)(4.76,-3.0092187)
\curveto(4.78,-3.0442188)(4.805,-3.1242187)(4.81,-3.1692188)
\curveto(4.815,-3.2142189)(4.825,-3.2992187)(4.83,-3.3392189)
\curveto(4.835,-3.3792188)(4.845,-3.4442186)(4.85,-3.4692187)
\curveto(4.855,-3.4942188)(4.86,-3.5542188)(4.86,-3.5892189)
\curveto(4.86,-3.6242187)(4.885,-3.7342188)(4.91,-3.8092186)
\curveto(4.935,-3.8842187)(4.975,-3.9992187)(4.99,-4.039219)
\curveto(5.005,-4.079219)(5.02,-4.144219)(5.02,-4.1692185)
\curveto(5.02,-4.1942186)(5.02,-4.244219)(5.02,-4.269219)
\curveto(5.02,-4.2942185)(5.02,-4.3242188)(5.02,-4.3392186)
}
\pscustom[linewidth=0.04,linestyle=dotted,dotsep=0.16cm]
{
\newpath
\moveto(3.64,-0.37921876)
\lineto(3.71,-0.32921875)
\curveto(3.745,-0.30421874)(3.825,-0.22421876)(3.87,-0.16921875)
\curveto(3.915,-0.11421875)(3.975,-0.02921875)(3.99,0.0)
\curveto(4.005,0.03078125)(4.055,0.12078125)(4.09,0.18078125)
\curveto(4.125,0.24078125)(4.19,0.39078125)(4.22,0.48078126)
\curveto(4.25,0.57078123)(4.315,0.78578126)(4.35,0.91078126)
\curveto(4.385,1.0357813)(4.435,1.3057812)(4.45,1.4507812)
\curveto(4.465,1.5957812)(4.49,1.8157812)(4.5,1.8907813)
\curveto(4.51,1.9657812)(4.53,2.1007812)(4.54,2.1607811)
\curveto(4.55,2.2207813)(4.56,2.3407812)(4.56,2.4007812)
\curveto(4.56,2.4607813)(4.565,2.5707812)(4.57,2.6207812)
\curveto(4.575,2.6707811)(4.585,2.7607813)(4.59,2.8007812)
\curveto(4.595,2.8407812)(4.605,2.9107811)(4.61,2.9407814)
\curveto(4.615,2.9707813)(4.62,3.0057812)(4.62,3.0207813)
}
\pscustom[linewidth=0.04,linestyle=dotted,dotsep=0.16cm]
{
\newpath
\moveto(4.1,0.12078125)
\lineto(4.13,0.16078125)
\curveto(4.145,0.18078125)(4.175,0.23078126)(4.19,0.26078126)
\curveto(4.205,0.29078126)(4.23,0.36078125)(4.24,0.40078124)
\curveto(4.25,0.44078124)(4.29,0.55078125)(4.32,0.62078124)
\curveto(4.35,0.69078124)(4.38,0.80078125)(4.38,0.8407813)
\curveto(4.38,0.88078123)(4.4,0.94578123)(4.42,0.97078127)
\curveto(4.44,0.99578124)(4.47,1.0657812)(4.48,1.1107812)
\curveto(4.49,1.1557813)(4.52,1.2407813)(4.54,1.2807813)
\curveto(4.56,1.3207812)(4.6,1.3957813)(4.62,1.4307812)
\curveto(4.64,1.4657812)(4.665,1.5857812)(4.67,1.6707813)
\curveto(4.675,1.7557813)(4.685,1.9057813)(4.69,1.9707812)
\curveto(4.695,2.0357811)(4.705,2.1407812)(4.71,2.1807814)
\curveto(4.715,2.2207813)(4.72,2.3257813)(4.72,2.3907812)
\curveto(4.72,2.4557812)(4.725,2.5807812)(4.73,2.6407812)
\curveto(4.735,2.7007813)(4.745,2.7957811)(4.75,2.8307812)
\curveto(4.755,2.8657813)(4.76,2.9307814)(4.76,2.9607813)
}
\pscustom[linewidth=0.04,linestyle=dotted,dotsep=0.16cm]
{
\newpath
\moveto(4.06,-0.35921875)
\lineto(4.11,-0.31921875)
\curveto(4.135,-0.29921874)(4.2,-0.24921875)(4.24,-0.21921875)
\curveto(4.28,-0.18921874)(4.34,-0.13921875)(4.36,-0.11921875)
\curveto(4.38,-0.09921875)(4.415,-0.04921875)(4.43,-0.01921875)
\curveto(4.445,0.01078125)(4.48,0.08078125)(4.5,0.12078125)
\curveto(4.52,0.16078125)(4.55,0.23078126)(4.56,0.26078126)
\curveto(4.57,0.29078126)(4.605,0.35078126)(4.63,0.38078126)
\curveto(4.655,0.41078126)(4.69,0.47578126)(4.7,0.5107812)
\curveto(4.71,0.54578125)(4.745,0.63078123)(4.77,0.68078125)
\curveto(4.795,0.73078126)(4.835,0.80078125)(4.85,0.82078123)
\curveto(4.865,0.8407813)(4.895,0.8957813)(4.91,0.93078125)
\curveto(4.925,0.9657813)(4.95,1.0307813)(4.96,1.0607812)
\curveto(4.97,1.0907812)(4.995,1.1707813)(5.01,1.2207812)
\curveto(5.025,1.2707813)(5.045,1.3657813)(5.05,1.4107813)
\curveto(5.055,1.4557812)(5.065,1.5957812)(5.07,1.6907812)
\curveto(5.075,1.7857813)(5.09,1.9957813)(5.1,2.1107812)
\curveto(5.11,2.2257812)(5.12,2.4057813)(5.12,2.4707813)
\curveto(5.12,2.5357811)(5.12,2.6707811)(5.12,2.7407813)
\curveto(5.12,2.8107812)(5.13,2.9057813)(5.14,2.9307814)
\curveto(5.15,2.9557812)(5.16,2.9857812)(5.16,3.0007813)
}
\pscustom[linewidth=0.04,linestyle=dotted,dotsep=0.16cm]
{
\newpath
\moveto(5.32,2.9007812)
\lineto(5.3,2.8107812)
\curveto(5.29,2.7657812)(5.28,2.6907814)(5.28,2.6607811)
\curveto(5.28,2.6307812)(5.28,2.5707812)(5.28,2.5407813)
\curveto(5.28,2.5107813)(5.28,2.4457812)(5.28,2.4107811)
\curveto(5.28,2.3757813)(5.275,2.3107812)(5.27,2.2807813)
\curveto(5.265,2.2507813)(5.25,2.1907814)(5.24,2.1607811)
\curveto(5.23,2.1307812)(5.215,2.0757813)(5.21,2.0507812)
\curveto(5.205,2.0257812)(5.19,1.9607812)(5.18,1.9207813)
\curveto(5.17,1.8807813)(5.155,1.8107812)(5.15,1.7807813)
\curveto(5.145,1.7507813)(5.14,1.6907812)(5.14,1.6607813)
\curveto(5.14,1.6307813)(5.14,1.5757812)(5.14,1.5507812)
\curveto(5.14,1.5257813)(5.14,1.4757812)(5.14,1.4507812)
\curveto(5.14,1.4257812)(5.14,1.3757813)(5.14,1.3507812)
\curveto(5.14,1.3257812)(5.135,1.2807813)(5.13,1.2607813)
\curveto(5.125,1.2407813)(5.105,1.1857812)(5.09,1.1507813)
\curveto(5.075,1.1157813)(5.045,1.0507812)(5.03,1.0207813)
\curveto(5.015,0.99078125)(4.99,0.93078125)(4.98,0.9007813)
\curveto(4.97,0.87078124)(4.94,0.82078123)(4.92,0.80078125)
\curveto(4.9,0.78078127)(4.87,0.74078125)(4.86,0.72078127)
\curveto(4.85,0.7007812)(4.835,0.65578127)(4.83,0.63078123)
\curveto(4.825,0.60578126)(4.815,0.55578125)(4.81,0.53078127)
\curveto(4.805,0.50578123)(4.785,0.44578126)(4.77,0.41078126)
\curveto(4.755,0.37578124)(4.73,0.31078124)(4.72,0.28078124)
\curveto(4.71,0.25078124)(4.68,0.19578125)(4.66,0.17078125)
\curveto(4.64,0.14578125)(4.61,0.09578125)(4.6,0.07078125)
\curveto(4.59,0.04578125)(4.565,0.00578125)(4.55,-0.00921875)
\curveto(4.535,-0.02421875)(4.51,-0.05921875)(4.5,-0.07921875)
\curveto(4.49,-0.09921875)(4.465,-0.12921876)(4.45,-0.13921875)
\curveto(4.435,-0.14921875)(4.4,-0.17421874)(4.38,-0.18921874)
\curveto(4.36,-0.20421875)(4.34,-0.22921875)(4.34,-0.25921875)
}
\pscustom[linewidth=0.04,linestyle=dotted,dotsep=0.16cm]
{
\newpath
\moveto(4.58,1.1407813)
\lineto(4.62,1.1907812)
\curveto(4.64,1.2157812)(4.675,1.2657813)(4.69,1.2907813)
\curveto(4.705,1.3157812)(4.72,1.3907813)(4.72,1.4407812)
\curveto(4.72,1.4907813)(4.73,1.5657812)(4.74,1.5907812)
\curveto(4.75,1.6157813)(4.76,1.6657813)(4.76,1.6907812)
\curveto(4.76,1.7157812)(4.765,1.7607813)(4.78,1.8207812)
}
\pscustom[linewidth=0.04,linestyle=dotted,dotsep=0.16cm]
{
\newpath
\moveto(3.66,-0.77921873)
\lineto(3.66,-0.87921876)
\curveto(3.66,-0.92921877)(3.66,-1.0292188)(3.66,-1.0792187)
\curveto(3.66,-1.1292187)(3.675,-1.2192187)(3.69,-1.2592187)
\curveto(3.705,-1.2992188)(3.73,-1.3642187)(3.74,-1.3892188)
\curveto(3.75,-1.4142188)(3.77,-1.5092187)(3.78,-1.5792187)
\curveto(3.79,-1.6492188)(3.81,-1.7642188)(3.82,-1.8092188)
\curveto(3.83,-1.8542187)(3.855,-1.9742187)(3.87,-2.0492187)
\curveto(3.885,-2.1242187)(3.91,-2.2842188)(3.92,-2.3692188)
\curveto(3.93,-2.4542189)(3.955,-2.5742188)(3.97,-2.6092188)
\curveto(3.985,-2.6442187)(4.005,-2.7942188)(4.01,-2.9092188)
\curveto(4.015,-3.0242188)(4.02,-3.1742187)(4.02,-3.2092187)
\curveto(4.02,-3.2442188)(4.02,-3.3892188)(4.02,-3.4992187)
\curveto(4.02,-3.6092188)(4.025,-3.7492187)(4.03,-3.7792187)
\curveto(4.035,-3.8092186)(4.035,-3.8792188)(4.03,-3.9192188)
\curveto(4.025,-3.9592187)(4.015,-4.034219)(4.01,-4.0692186)
\curveto(4.005,-4.104219)(4.0,-4.159219)(4.0,-4.2192187)
}
\pscustom[linewidth=0.04,linestyle=dotted,dotsep=0.16cm]
{
\newpath
\moveto(3.76,-0.9792187)
\lineto(3.79,-1.1792188)
\curveto(3.805,-1.2792188)(3.84,-1.4342188)(3.86,-1.4892187)
\curveto(3.88,-1.5442188)(3.91,-1.6392188)(3.92,-1.6792188)
\curveto(3.93,-1.7192187)(3.95,-1.8192188)(3.96,-1.8792187)
\curveto(3.97,-1.9392188)(3.99,-2.0542188)(4.0,-2.1092188)
\curveto(4.01,-2.1642187)(4.03,-2.2592187)(4.04,-2.2992187)
\curveto(4.05,-2.3392189)(4.065,-2.4092188)(4.07,-2.4392188)
\curveto(4.075,-2.4692187)(4.085,-2.5492187)(4.09,-2.5992188)
\curveto(4.095,-2.6492188)(4.1,-2.7342188)(4.1,-2.7692187)
\curveto(4.1,-2.8042188)(4.115,-2.8992188)(4.13,-2.9592187)
\curveto(4.145,-3.0192187)(4.165,-3.1142187)(4.17,-3.1492188)
\curveto(4.175,-3.1842186)(4.185,-3.2492187)(4.19,-3.2792187)
\curveto(4.195,-3.3092186)(4.2,-3.3992188)(4.2,-3.4592187)
\curveto(4.2,-3.5192187)(4.195,-3.5992188)(4.19,-3.6192188)
\curveto(4.185,-3.6392188)(4.18,-3.6992188)(4.18,-3.7392187)
\curveto(4.18,-3.7792187)(4.175,-3.8442187)(4.17,-3.8692188)
\curveto(4.165,-3.8942187)(4.155,-3.9442186)(4.15,-3.9692187)
\curveto(4.145,-3.9942188)(4.14,-4.0442185)(4.14,-4.0692186)
\curveto(4.14,-4.0942187)(4.14,-4.144219)(4.14,-4.1692185)
}
\pscustom[linewidth=0.04,linestyle=dotted,dotsep=0.16cm]
{
\newpath
\moveto(3.82,-0.67921877)
\lineto(3.85,-0.70921874)
\curveto(3.865,-0.7242187)(3.89,-0.76421875)(3.9,-0.7892187)
\curveto(3.91,-0.81421876)(3.94,-0.8542187)(3.96,-0.86921877)
\curveto(3.98,-0.88421875)(4.005,-0.92921877)(4.01,-0.95921874)
\curveto(4.015,-0.9892188)(4.035,-1.0542188)(4.05,-1.0892187)
\curveto(4.065,-1.1242187)(4.12,-1.2342187)(4.16,-1.3092188)
\curveto(4.2,-1.3842187)(4.245,-1.4792187)(4.25,-1.4992187)
\curveto(4.255,-1.5192188)(4.27,-1.5742188)(4.28,-1.6092187)
\curveto(4.29,-1.6442188)(4.32,-1.7242187)(4.34,-1.7692188)
\curveto(4.36,-1.8142188)(4.4,-1.9342188)(4.42,-2.0092187)
\curveto(4.44,-2.0842187)(4.465,-2.2042189)(4.47,-2.2492187)
\curveto(4.475,-2.2942188)(4.485,-2.3842187)(4.49,-2.4292188)
\curveto(4.495,-2.4742188)(4.505,-2.5742188)(4.51,-2.6292188)
\curveto(4.515,-2.6842186)(4.53,-2.8142188)(4.54,-2.8892188)
\curveto(4.55,-2.9642189)(4.565,-3.1142187)(4.57,-3.1892188)
\curveto(4.575,-3.2642188)(4.58,-3.3742187)(4.58,-3.4092188)
\curveto(4.58,-3.4442186)(4.57,-3.5242188)(4.56,-3.5692186)
\curveto(4.55,-3.6142187)(4.53,-3.6942186)(4.52,-3.7292187)
\curveto(4.51,-3.7642188)(4.5,-3.8242188)(4.5,-3.8492188)
\curveto(4.5,-3.8742187)(4.495,-3.9242187)(4.49,-3.9492188)
\curveto(4.485,-3.9742188)(4.475,-4.0242186)(4.47,-4.0492187)
\curveto(4.465,-4.0742188)(4.46,-4.124219)(4.46,-4.1492186)
}
\pscustom[linewidth=0.04,linestyle=dotted,dotsep=0.16cm]
{
\newpath
\moveto(4.62,-2.7192187)
\lineto(4.62,-2.8092186)
\curveto(4.62,-2.8542187)(4.63,-2.9392188)(4.64,-2.9792187)
\curveto(4.65,-3.0192187)(4.665,-3.0892189)(4.67,-3.1192188)
\curveto(4.675,-3.1492188)(4.685,-3.2342188)(4.69,-3.2892187)
\curveto(4.695,-3.3442187)(4.705,-3.4692187)(4.71,-3.5392187)
\curveto(4.715,-3.6092188)(4.73,-3.7192187)(4.74,-3.7592187)
\curveto(4.75,-3.7992187)(4.76,-3.9092188)(4.76,-3.9792187)
\curveto(4.76,-4.0492187)(4.76,-4.144219)(4.76,-4.2192187)
}
\pscustom[linewidth=0.04,linestyle=dotted,dotsep=0.16cm]
{
\newpath
\moveto(4.62,-3.6792188)
\lineto(4.61,-3.7792187)
\curveto(4.605,-3.8292189)(4.6,-3.9142187)(4.6,-3.9492188)
\curveto(4.6,-3.9842188)(4.6,-4.0492187)(4.6,-4.079219)
}
\pscustom[linewidth=0.04,linestyle=dotted,dotsep=0.16cm]
{
\newpath
\moveto(4.46,-3.2792187)
\lineto(4.44,-3.3592188)
\curveto(4.43,-3.3992188)(4.415,-3.4792187)(4.41,-3.5192187)
\curveto(4.405,-3.5592186)(4.4,-3.6242187)(4.4,-3.6492188)
\curveto(4.4,-3.6742187)(4.4,-3.7342188)(4.4,-3.7692187)
\curveto(4.4,-3.8042188)(4.395,-3.8692188)(4.39,-3.8992188)
\curveto(4.385,-3.9292188)(4.38,-3.9792187)(4.38,-4.039219)
}
\pscustom[linewidth=0.04,linestyle=dotted,dotsep=0.16cm]
{
\newpath
\moveto(4.02,-1.6192187)
\lineto(4.02,-1.7292187)
\curveto(4.02,-1.7842188)(4.025,-1.8842187)(4.03,-1.9292188)
\curveto(4.035,-1.9742187)(4.05,-2.0542188)(4.06,-2.0892189)
\curveto(4.07,-2.1242187)(4.085,-2.1792188)(4.09,-2.1992188)
\curveto(4.095,-2.2192187)(4.11,-2.2642188)(4.12,-2.2892187)
\curveto(4.13,-2.3142188)(4.15,-2.3742187)(4.16,-2.4092188)
\curveto(4.17,-2.4442186)(4.185,-2.5142188)(4.19,-2.5492187)
\curveto(4.195,-2.5842187)(4.205,-2.6542187)(4.21,-2.6892188)
\curveto(4.215,-2.7242188)(4.22,-2.7842188)(4.22,-2.8092186)
\curveto(4.22,-2.8342187)(4.22,-2.9142187)(4.22,-2.9692187)
\curveto(4.22,-3.0242188)(4.22,-3.0842187)(4.22,-3.0992188)
}
\pscustom[linewidth=0.04,linestyle=dotted,dotsep=0.16cm]
{
\newpath
\moveto(3.92,-2.5992188)
\lineto(3.92,-2.6792188)
\curveto(3.92,-2.7192187)(3.925,-2.7892187)(3.93,-2.8192186)
\curveto(3.935,-2.8492188)(3.945,-2.9042187)(3.95,-2.9292188)
\curveto(3.955,-2.9542189)(3.96,-3.0242188)(3.96,-3.0692186)
\curveto(3.96,-3.1142187)(3.96,-3.1842186)(3.96,-3.2092187)
\curveto(3.96,-3.2342188)(3.955,-3.2842188)(3.95,-3.3092186)
\curveto(3.945,-3.3342187)(3.94,-3.3842187)(3.94,-3.4092188)
\curveto(3.94,-3.4342186)(3.94,-3.4842188)(3.94,-3.5092187)
\curveto(3.94,-3.5342188)(3.945,-3.5992188)(3.95,-3.6392188)
\curveto(3.955,-3.6792188)(3.96,-3.7492187)(3.96,-3.7792187)
\curveto(3.96,-3.8092186)(3.96,-3.8642187)(3.96,-3.8892188)
\curveto(3.96,-3.9142187)(3.96,-3.9642189)(3.96,-3.9892187)
\curveto(3.96,-4.014219)(3.955,-4.059219)(3.95,-4.079219)
}
\pscustom[linewidth=0.04,linestyle=dotted,dotsep=0.16cm]
{
\newpath
\moveto(3.62,-0.75921875)
\lineto(3.59,-0.88921875)
\curveto(3.575,-0.95421875)(3.56,-1.0742188)(3.56,-1.1292187)
\curveto(3.56,-1.1842188)(3.565,-1.2692188)(3.57,-1.2992188)
\curveto(3.575,-1.3292187)(3.595,-1.3592187)(3.61,-1.3592187)
\curveto(3.625,-1.3592187)(3.655,-1.3842187)(3.67,-1.4092188)
\curveto(3.685,-1.4342188)(3.715,-1.4792187)(3.73,-1.4992187)
\curveto(3.745,-1.5192188)(3.76,-1.5442188)(3.76,-1.5592188)
}
\pscustom[linewidth=0.04,linestyle=dotted,dotsep=0.16cm]
{
\newpath
\moveto(3.96,-1.0192188)
\lineto(3.97,-1.0792187)
\curveto(3.975,-1.1092187)(3.99,-1.1642188)(4.0,-1.1892188)
\curveto(4.01,-1.2142187)(4.03,-1.2642188)(4.06,-1.3392187)
}
\pscustom[linewidth=0.04,linestyle=dotted,dotsep=0.16cm]
{
\newpath
\moveto(3.92,0.06078125)
\lineto(3.96,0.12078125)
\curveto(3.98,0.15078124)(4.02,0.20578125)(4.04,0.23078126)
\curveto(4.06,0.25578126)(4.09,0.31578124)(4.1,0.35078126)
\curveto(4.11,0.38578126)(4.125,0.45578125)(4.13,0.49078125)
\curveto(4.135,0.5257813)(4.15,0.59578127)(4.16,0.63078123)
\curveto(4.17,0.66578126)(4.185,0.73078126)(4.19,0.7607812)
\curveto(4.195,0.79078126)(4.215,0.87578124)(4.23,0.93078125)
\curveto(4.245,0.98578125)(4.275,1.0607812)(4.29,1.0807812)
\curveto(4.305,1.1007812)(4.32,1.1307813)(4.32,1.1607813)
}
\pscustom[linewidth=0.04,linestyle=dotted,dotsep=0.16cm]
{
\newpath
\moveto(4.28,-1.7192187)
\lineto(4.28,-1.8192188)
\curveto(4.28,-1.8692187)(4.28,-1.9592187)(4.28,-1.9992187)
\curveto(4.28,-2.0392187)(4.29,-2.1042187)(4.3,-2.1292188)
\curveto(4.31,-2.1542187)(4.335,-2.2092187)(4.35,-2.2392187)
\curveto(4.365,-2.2692187)(4.385,-2.3292189)(4.39,-2.3592188)
\curveto(4.395,-2.3892188)(4.415,-2.4492188)(4.43,-2.4792187)
\curveto(4.445,-2.5092187)(4.465,-2.5542188)(4.48,-2.5992188)
}
\pscustom[linewidth=0.04,linestyle=dotted,dotsep=0.16cm]
{
\newpath
\moveto(4.16,-1.2592187)
\lineto(4.17,-1.3192188)
\curveto(4.175,-1.3492187)(4.2,-1.4442188)(4.22,-1.5092187)
\curveto(4.24,-1.5742188)(4.28,-1.6792188)(4.3,-1.7192187)
\curveto(4.32,-1.7592187)(4.36,-1.8192188)(4.38,-1.8392187)
\curveto(4.4,-1.8592187)(4.42,-1.9192188)(4.42,-1.9592187)
\curveto(4.42,-1.9992187)(4.425,-2.0642188)(4.43,-2.0892189)
\curveto(4.435,-2.1142187)(4.45,-2.1742187)(4.46,-2.2092187)
\curveto(4.47,-2.2442188)(4.485,-2.2992187)(4.5,-2.3592188)
}
\pscustom[linewidth=0.04,linestyle=dotted,dotsep=0.16cm]
{
\newpath
\moveto(4.52,1.3807813)
\lineto(4.53,1.6207813)
\curveto(4.535,1.7407813)(4.54,1.9307812)(4.54,2.0007813)
\curveto(4.54,2.0707812)(4.54,2.1957812)(4.54,2.2507813)
\curveto(4.54,2.3057814)(4.565,2.3957813)(4.59,2.4307814)
\curveto(4.615,2.4657812)(4.65,2.5857813)(4.66,2.6707811)
\curveto(4.67,2.7557812)(4.685,2.8607812)(4.7,2.9207811)
}
\pscustom[linewidth=0.04,linestyle=dotted,dotsep=0.16cm]
{
\newpath
\moveto(3.72,-0.93921876)
\lineto(3.73,-1.0192188)
\curveto(3.735,-1.0592188)(3.765,-1.1392188)(3.79,-1.1792188)
\curveto(3.815,-1.2192187)(3.85,-1.2942188)(3.86,-1.3292187)
\curveto(3.87,-1.3642187)(3.89,-1.4442188)(3.9,-1.4892187)
\curveto(3.91,-1.5342188)(3.93,-1.6392188)(3.94,-1.6992188)
\curveto(3.95,-1.7592187)(3.975,-1.8942188)(3.99,-1.9692187)
\curveto(4.005,-2.0442188)(4.025,-2.1292188)(4.04,-2.1592188)
}
\pscustom[linewidth=0.04,linestyle=dotted,dotsep=0.16cm]
{
\newpath
\moveto(4.06,-2.6392188)
\lineto(4.04,-2.8292189)
\curveto(4.03,-2.9242187)(4.02,-3.0892189)(4.02,-3.1592188)
\curveto(4.02,-3.2292187)(4.02,-3.3542187)(4.02,-3.4092188)
\curveto(4.02,-3.4642189)(4.03,-3.5492187)(4.04,-3.5792189)
\curveto(4.05,-3.6092188)(4.07,-3.6692188)(4.08,-3.6992188)
}
\pscustom[linewidth=0.04,linestyle=dotted,dotsep=0.16cm]
{
\newpath
\moveto(4.58,-3.2392187)
\lineto(4.59,-3.3292189)
\curveto(4.595,-3.3742187)(4.635,-3.4692187)(4.67,-3.5192187)
\curveto(4.705,-3.5692186)(4.755,-3.6492188)(4.77,-3.6792188)
\curveto(4.785,-3.7092187)(4.805,-3.7792187)(4.81,-3.8192186)
\curveto(4.815,-3.8592188)(4.825,-3.9292188)(4.83,-3.9592187)
}
\pscustom[linewidth=0.04,linestyle=dotted,dotsep=0.16cm]
{
\newpath
\moveto(3.9,-1.7192187)
\lineto(3.92,-1.9092188)
\curveto(3.93,-2.0042188)(3.955,-2.1742187)(3.97,-2.2492187)
\curveto(3.985,-2.3242188)(4.005,-2.4492188)(4.01,-2.4992187)
\curveto(4.015,-2.5492187)(4.025,-2.6142187)(4.04,-2.6592188)
}
\psline[linewidth=0.04cm,arrowsize=0.05291667cm 2.0,arrowlength=1.4,arrowinset=0.4]{<-}(5.12,3.1407812)(5.98,4.0607815)
\psline[linewidth=0.04cm,linecolor=color282,linestyle=dashed,dash=0.16cm 0.16cm,arrowsize=0.05291667cm 2.0,arrowlength=1.4,arrowinset=0.4]{<-}(2.3,0.22078125)(1.84,3.7407813)
\psline[linewidth=0.04cm,linecolor=color282,linestyle=dashed,dash=0.16cm 0.16cm,arrowsize=0.05291667cm 2.0,arrowlength=1.4,arrowinset=0.4]{<-}(4.88,2.2607813)(1.9,3.6407812)
\psline[linewidth=0.04cm,linestyle=dashed,dash=0.16cm 0.16cm,arrowsize=0.05291667cm 2.0,arrowlength=1.4,arrowinset=0.4]{<-}(4.92,0.8407813)(7.92,-0.77921873)
\psline[linewidth=0.04cm,linestyle=dashed,dash=0.16cm 0.16cm,arrowsize=0.05291667cm 2.0,arrowlength=1.4,arrowinset=0.4]{<-}(4.66,-3.6192188)(7.96,-0.89921874)
\usefont{T1}{ptm}{m}{n}
\rput(6.36625,4.210781){V}
\usefont{T1}{ptm}{m}{n}
\rput(8.481406,-0.86921877){$\cN$}
\usefont{T1}{ptm}{m}{n}
\rput(1.7295313,3.9307814){\color{blue}M($\psi $)}
\usefont{T1}{ptm}{m}{n}
\rput(3.4065626,3.2307813){\color{color282}A}
\usefont{T1}{ptm}{m}{n}
\rput(1.7671875,1.9907813){\color{color420}B}
\end{pspicture} 
}
\end{center}
\vspace{0.3cm}

\begin{remark}\label{r:pichonseade}
In \cite[Theorem 1.3]{PS} the authors state that if $\Sing \psi \subset V$ and if Thom's condition holds along V, then there is an empty tube fibration like (\ref{eq:emptytube}) which blows out to a Milnor fibration on the sphere. As the authors remark themselves \cite[p. 488 and 492]{PS},  this is a re-formulation of Milnor's results which follows from Milnor's method. There are the following differences with our Theorem \ref{t:milnor}: 
\begin{itemize}
\item[(1).]  \sloppy the Thom regularity hypothesis of \cite[Theorem 1.3]{PS} implies our condition $\overline{M(\psi)\m V} \cap  V = \{ 0\}$ whereas the converse is presumably not true in general (see below \S \ref{ss:thom} for a discussion of this relation and Example \ref{e:thom}).
\item[(2).] \cite[Theorem 1.3]{PS} exploits only the fibration outside a tube. We here need  the more refined fibration structure inside the tube, which is   necessary for having an open book structure as in Definition \ref{d:booksing}.  
\end{itemize}
\end{remark}


\section{Proof of Theorem \ref{t:isolcrt} and some consequences} \label{s:strong}

We assume here condition \eqref{eq:main} but we do not assume $\Sing \psi \subset V$.
As explained above at \S \ref{ss:endproof}, the condition  $\overline{M(\psi)\m V} \cap  V = \{ 0\}$ is equivalent to the existence of a conical neighbourhood $\cN$ of $V \m \{ 0\}$ such that $M(\psi) \cap \cN \m V = \emptyset$. This implies neither the existence of a ``full tube'' fibration (\ref{eq:tube}),  nor the existence of the ``empty tube fibration'' (\ref{eq:emptytube}), since our hypothesis \eqref{eq:main} allows singularities on the fibres of $\psi$ close to $V$. It is nevertheless enough to imply that the restriction:
\begin{equation}\label{eq:tube-sphere}
\psi_| : S^{m-1}_\e \cap \psi^{-1}( \bar B^{p}_\eta \m \{0\}) \to   \bar B^{p}_\eta \m \{0\}
\end{equation}
is a proper submersion, for all small enough $\e \gg \eta>0$, hence surjective too.

On the other hand,
the condition that the germ of $M(\frac{\psi}{\|\psi\|})$ at the origin is empty means that the mapping 
$\frac{\psi}{\| \psi\|} : S^{m-1}_\e \m \psi^{-1}( B^{p}_\eta) \to  S^{p-1}_1$ is a proper submersion for any small enough $\e >0$. Following the proof of Theorem \ref{t:milnor}, we notice that we have here the particularity that the mappings to be glued along the boundary $S^{m-1}_\e \cap \psi^{-1}( S^{p-1}_\eta)$ are induced by the same mapping $\frac{\psi}{\|\psi\|}$ from both sides, hence this glueing is trivial, and of course smooth, with $\theta := \frac{\psi}{\|\psi\|}$.
 The conclusion follows now just like in the proof of Theorem \ref{t:milnor}.
\fin

\begin{example}\label{r:ps}
It was shown in \cite{PS} that the mappings of type  $f\bar g : (\bC^2, 0) \to (\bC, 0)$, where $f, g :\bC^2 \to \bC$ are holomorphic functions such that $\Sing f\bar g \subset V$ have a Milnor fibration  $f\bar g/ | f\bar g| : S^3_\e \m K_\e \to S^1$ and are Thom regular along $V$.  Since Thom regularity implies condition \eqref{eq:main} (see also \S \ref{ss:thom}), these mean that the hypotheses of our Theorem \ref{t:isolcrt} are satisfied, hence one obtains open book structures with singular binding $(K_\e, f\bar g/ | f\bar g|)$.
\end{example}



\subsection{Radial weighted-homogeneous real mappings}\label{ss:radial}
Let us consider the $\bR_+$-action on $\bR^m$:
$\rho\cdot x  =(\rho^{q_{1}}x_{1},\ldots,\rho^{q_{m}}x_{m})$ for $\rho \in \bR_+$ and $q_1, \ldots, q_m \in \bN^*$ relatively prime positive integers. Let $\gamma (x) := \sum_{j=1}^m q_j x_j \frac{\partial}{\partial x_j}$ be the corresponding Euler vector field on $\bR^m$; we have $\gamma (x) = 0$ if and only if $x=0$.

We say that the mapping $\psi$ is \textit{radial weighted-homogeneous} of degree $d >0$ if $\psi (\rho \cdot x) = \rho^d\psi(x)$ for all $x$ in some neighbourhood of $0$.

\begin{proposition} \label{p:radial}
  If $\psi$ is radial weighted-homogeneous and $\Sing \psi \subset V$, then $M(\frac{\psi}{\|\psi\|}) = \emptyset$. 
\end{proposition}
\begin{proof} Let us first remark that $\Sing \psi \subset V$ implies that $\frac{\psi}{\|\psi\|}: B^m_{\e}\setminus V \to S^{p-1}_1$ is a submersion.
We use the following criterion equivalent to $\{ 0\} \not\in M(\frac{\psi}{\|\psi\|})$ from the proof of \cite[Theorem 2.2]{dST1}:
 $\exists \e_0 >0$  such that  $\rank  \Omega_\psi(x)  = p$,  $\forall x\in B^m_{\e_0}\setminus V$,
where  $\Omega_\psi(x)$ denotes the $[(p-1)p/2 +1] \times m$ matrix defined as: 
\[ \Omega_\psi(x) := \left[ \begin{array}{cc}
         \omega_{1,2}(x) \\ \vdots \\ \omega_{i,j}(x) \\ \vdots  \\ \omega_{p-1, p}(x) \\
                       x_1, \ldots , x_m  
 \end{array} \right]    \]
 having in each of the rows the vector $\omega_{i,j}(x) := \psi_i(x)\hspace{1pt} \grad \psi_j(x) - \psi_j(x)\hspace{1pt} \grad \psi_i(x)$, for $i,j = 1, \ldots ,p$ with $i<j$, except for the last row which contains the position vector  $(x_1, \ldots , x_m)$. 
Observing that $\langle \gamma(x),  \grad \psi_i(x) \rangle = d \cdot \psi_i(x)$ for any $i$ and any $x\in B^m_{\e_0}\setminus V$, we have:
\[ \langle \gamma(x),  \omega_{ij}(x) \rangle = d [\psi_i(x)\psi_j(x) - \psi_j(x)\psi_i(x)] = 0,
 \]
which means that the Euler vector field $\gamma(x)$ is tangent to the fibres of $\frac{\psi}{\|\psi\|}$.
We also have:
 \begin{equation}\label{eq:ineg}
   \langle \gamma(x),  x \rangle = \sum_i q_i x_i^2 > 0,
 \end{equation}
for $x \not= 0$. This shows that the position vector $x$ cannot be orthogonal to the tangent space of the fibres of $\frac{\psi}{\|\psi\|}$, which means that the sphere $S^{m-1}_\e$ is transverse to the fibres of $\frac{\psi}{\|\psi\|}$, for any 
$\e >0$.  
\end{proof}

We get the following statement, the proof of which is an immediate consequence of the proof of Theorem \ref{t:isolcrt} via Proposition \ref{p:radial}. It also represents an extension to nonisolated singularities of our previous \cite[Theorem 3.1]{dST1}, see also the proof of \cite[Theorem 4.1]{dST0} and compare with \cite[Example 2.1.4]{CSS}.

\begin{corollary}\label{c:radial}
  Let $\psi$ be radial weighted-homogeneous with $\Sing \psi \subset V$ and satisfying condition \eqref{eq:main} with $\codim V = p$. Then  $(K_\e, \frac{\psi}{\|\psi \|})$ is a (singular) open book
 decomposition.
\fin
\end{corollary}

\section{Polar weighted-homogeneous mixed functions}\label{ss:polar}

We consider a mixed polynomial\footnote{the term ``mixed polynomial'' was introduced by Oka \cite{oka1}.} $f(\mathbf{z},\mathbf{\overline{z}})=\sum_{v,\mu}c_{\nu,\mu}\mathbf{z}^{\nu}\overline{\mathbf{z}}^{\mu}$ where $\mathbf{z}=(z_{1},\ldots,z_{n})$, $\mathbf{\overline{z}}=(\overline{z}_{1},\ldots,\overline{z}_{n})$,
$\mathbf{z}^{\nu}=z_{1}^{\nu_{1}}\cdots z_{n}^{\nu_{n}}$, $\overline{\mathbf{z}}^{\mu}=\overline{z}_{1}^{\mu_{1}}\cdots\overline{z}_{n}^{\mu_{n}}$
for $\nu=(\nu_{1},\ldots,\nu_{n})$ and $\mu=(\mu_{1},\ldots,\mu_{n})$ non-negative integer exponents.
 As a matter of fact, any mixed polynomial is a real polynomial mapping $\bR^{2n} \to \bR^2$, and conversely.

After \cite{CM} and \cite{oka1},  $f$ is called \emph{polar weighted-homogeneous} if there are non-zero integers 
$p_{1},\ldots,p_{n}$ and $k$ such that $\gcd(p_{1},\ldots p_{n})=1$ and $\sum_{j=1}^{n}p_{j}(\nu_{j}-\mu_{j})=k$. The corresponding  $S^{1}$-action on $\mathbb{C}^{n}$ is: 
\[ 
\lambda\cdot(\mathbf{z},\mathbf{\overline{z}})  =(\lambda^{p_{1}}z_{1},\ldots,\lambda^{p_{n}}z_{n},\lambda^{-p_{1}}\overline{z}_{1},\ldots,\lambda^{-p_{n}}\overline{z}_{n}),\,\lambda\in S^{1}.
\]
Notice that ``polar weighted-homogeneous'' and ``radial weighted-homogeneous''are two independent notions.

\subsection{Proof of Theorem \ref{t:homogen}}
We first prove that $\im f$ contains a small enough disk at $0 \in \bC$. Since  $f\not\equiv 0$ and $\im f$ is a semi-algebraic set germ, by the Curve Selection Lemma, the image contains a curve $\gamma$ which intersects the circles $S^1_\eta \subset \bC$ for any small enough radius $\eta>0$. 
Take now some $a\in S^1_\eta \cap \gamma$ and $z\in f^{-1}(a)$. Since $f(\lambda \cdot (\mathbf{z}, \mathbf{\overline{z}})) = \lambda^k f(\mathbf{z}, \mathbf{\overline{z}})$, we have $\lambda^k a \in \im f$ for any $\lambda \in S^1$. This shows that $\im f$ contains a disk $D^2_{\eta_0}$, for some small enough $\eta_0>0$.

The germ at $0$ of the set of critical values of $f$ is a semi-algebraic set of dimension $\le 1$. Take its complement in $\im f$, which is a 2 dimensional semi-algebraic germ at $0$. Applying the above reasoning yields that all values $\not= 0$ are regular, hence their inverse images are manifolds of dimension $2n-2\ge 1$.

Take now the restriction of $f$ to some small enough sphere $S_\e^{m-1}$. Its image must contain a non-constant curve germ at $0$. Since the $S^1$-action preserves the sphere, by the same reasoning as above, the image of $f_{|S_\e^{m-1}}$ contains a disk $D^2_{\eta_0}$. The regular values of $f_{|S_\e^{m-1}}$ are a dense semi-algebraic set and if $a$ is a regular value then $\lambda^k a$ is regular too, for any $\lambda \in S^1$. Hence all values of the pointed disk $D^2_{\eta_0}\m \{ 0\}$ are regular. This shows that we have a conical neighbourhood $\cN$ of $V$ such that the $\rho$-regularity holds within $\cN \m V$, hence the property \eqref{eq:main} holds.

In order to prove that $M(\frac{f}{\|f\|}) = \emptyset$ we apply the same reasoning to the mapping $\phi := \frac{f}{\|f\|}$ since the $S^1$-action yields $\phi(\lambda \cdot (\mathbf{z}, \mathbf{\overline{z}})) = \lambda^k \phi(\mathbf{z}, \mathbf{\overline{z}})$ and preserves the spheres centred at the origin.
 Our statement follows now from Theorem \ref{t:isolcrt}.
\fin

\begin{remark}\label{r:cisneros}
\sloppy In addition to our hypothesis of ``polar weighted-homogeneous'', if one assumes that  $f$ is moreover radial weighted-homogeneous and that all weights are positive, then the Milnor fibration induced by $f/|f|$ on the spheres was observed by Oka \cite[\S 5.4]{oka1} and Cisneros-Molina \cite[Prop. 3.4]{CM}. These yield property (a) from \S \ref{intro} in this particular setting. 

Note that our Theorem \ref{t:homogen} drops half of the hypotheses of the preceding results in what  concerns property (a) of \S \ref{intro}, and that our proof has a different flavor. Moreover, it addresses the issue (b) of \S \ref{intro} since, according to our definition, the open book structure with singularities in the binding follows only from the conjunction of properties (a) and (b) and not only from one of them.  See \S \ref{ss:thom} for remarks concerning the conditions which one has to impose near the link, Oka's result \cite[Theorem 52]{oka2} and Example \ref{e:thom}.
\end{remark}

\begin{example}\label{e:hom}
$f(x,y) = x^4\bar y^2 + y^2$ is not radial weighted homogeneous but a polar homogeneous mixed function germ $(\bC^2,0) \to (\bC,0)$, and we have $\Sing f = V = \{ y=0\}$.  According to Theorem \ref{t:homogen}, $f/ | f|$ defines  on $S^3$ an open book structure with singular binding. 
\end{example}
\section{Constructing more examples with nonisolated singularities}\label{s:classicandnew}

\subsection{Condition \eqref{eq:main} and Thom regularity condition}\label{ss:thom}

It is well-known that if $V = \psi^{-1}(0)$ may be endowed with a stratification $\cS$ such that in a sufficiently small ball $B_{\epsilon}$ the pair $(B_{\epsilon}\setminus V, S)$ satisfies Thom's condition (a$_\psi$), for all $S\in \cS $, then the stratified transversality of $V$ to 
all small enough spheres implies the transversality of the spheres to the nearby fibres in some neighbourhood $\cN$ of $V\m \{ 0\}$. This transversality is equivalent to the condition \eqref{eq:main} and by the proof of Theorem \ref{t:isolcrt} we conclude the existence of the full tube fibration on all small enough spheres \eqref{eq:tube-sphere}.
 This  well-known
observation may be traced back at least to \cite{HL}\footnote{see also Hironaka \cite{Hi} and \cite{Le-Oslo}.}.

It is conjecturally possible that the condition \eqref{eq:main} does not imply the Thom (a$_\psi$)-regularity, but we could not find an example yet. Let us at least point out two situations where one proves directly the weaker condition \eqref{eq:main}. One of them is contained in the proof of \cite[Lemma 51]{oka2} in the setting of mixed functions $f$, where Oka shows the existence of a full tube fibration for a special class of mappings, namely the ``super strongly non-degenerate mixed functions''. This is a condition which allows Oka to prove in \cite[Theorem 52]{oka2} the existence of the Milnor fibration on spheres (by extending Milnor's method)
and its equivalence to the tube fibration. Altogether these properties yield, in our terminology, an open book structure with singular binding induced by the mapping $\frac{f}{| f| }$.

Another example of computation is the following one, where the mapping is not a mixed function, hence this is different from all the previously mentioned situations of \cite[Theorem 52]{oka2} or of Theorem \ref{t:homogen}. 

\begin{example}\label{e:thom}
Let $f :\mathbb{R}^{3}\rightarrow\mathbb{R}^{2}$,  $f(x,y,z)=(y^{4}-z^{2}x^{2}-x^{4},xy)$. Then $V(f)$ is the real line $\{(x,y,z)\in\mathbb{R}^{3}\mid x=y=0\}$ and $\Sing f=V(f)$.

Let us show \eqref{eq:main}. If this were not true, then, by the Curve Selection Lemma, there are some analytic curves $x(t)$, $y(t)$, $z(t)$, $a(t)$ and $b(t)$ defined on a
small enough interval $\left]0,\varepsilon\right[$ such that $\lim_{t\to 0} x(t) =\lim_{t\to 0} y(t)=0$,
$\lim_{t\to 0} z(t)=z_{0}\neq 0$ and 
\[ \begin{array}{l}
x(t)  =  a(t)(-4x^{3}(t)-2x(t)z^{2}(t))+b(t)y(t),\\
y(t)  =  4a(t)y^{3}(t)+b(t)x(t),\\
z(t)  =  -2a(t)z(t)x^{2}(t).
\end{array} \]

Let  $x(t)=x_{0}t^{\beta}+\mathrm{h.o.t.}$, where $x_{0}\neq0$ and $\beta\in\mathbb{N}$. From the third line and from $\lim_{t\to 0} z(t) = z_{0}$
we get $a(t)=-\frac{1}{2x_{0}^{2}}t^{-2\beta}+\mathrm{h.o.t.}$
We eliminate $b(t)$ from the first two lines and get:
\[y^{2}(t)-x^{2}(t)=a(t)(4y^{4}(t)+4x^{4}(t)+2x^{2}(t)z^{2}(t)).\]

From this equality, since
 $\lim_{t\to 0}2 x^{2}(t)z^{2}(t)a(t)=-z_{0}^{2}<0$ and $\mathrm{ord_{t}}(a(t) x^{4}(t)) = 2 \beta >0$,
and $\lim_{t\to 0}(y^{2}(t)-x^{2}(t))=0$, we get 
$y^{4}(t)= - \frac{1}{2}  x_0^2 z_0^2 t^{2 \beta}+\mathrm{h.o.t.}$, which yields a sign contradiction. 

\end{example}

\subsection{A real Thom-Sebastiani type statement}\label{ss:ThomSebastiani}
 In order to build further examples, we show the following. 
\begin{proposition} \label{p:exempl}
 Consider two mappings in separate variables,  $\psi : (\bR^m,0) \to (\bR^p,0)$ and $\phi : (\bR^n,0) \to (\bR^p,0)$, such that $\Sing \psi \subset V(\psi)$ and $\Sing \phi \subset V(\phi)$, and both $V(\psi)$ and $V(\phi)$ have codimension $p$. Assume that $\psi$ and $\phi$ satisfy the Thom regularity condition at $V(\psi)$ and $V(\phi)$, respectively.

Then  $\psi + \phi : (\bR^m \times \bR^n, 0) \to (\bR^p,0)$ satisfies the Thom regularity condition and there exists a higher open book structure with singular binding $(K_{\psi +\phi}, \theta)$ on $S^{m+p-1}_\e$, which is
   independent of $\e >0$ small enough, up to isotopies.

If moreover 
 $\psi$ and $\phi$ are radial weighted-homogeneous  then $(K_{\psi + \phi}, \frac{\psi + \phi}{\| \psi + \phi \| })$ is a higher open book with singular binding.  

\end{proposition}
\begin{proof} 
From $\Sing \psi \subset V(\psi)$ and $\Sing \phi \subset V(\phi)$ it follows, by checking the rank of the Jacobian matrix, that $\Sing (\psi + \phi) \subset \Sing \psi \times \Sing \phi \subset V(\psi) \times  V(\phi) \subset  V(\psi + \phi) \subset \bR^m \times \bR^n$. 
The sum of separate variables mappings which both satisfy the Thom regularity condition has the same property. Indeed, if $\cW_1$ and $\cW_2$ denote some Thom stratifications of $V(\psi)$ and $V(\phi)$ respectively, then the product stratification $\cW_1 \times \cW_2$ satisfies the Thom ($a_{\psi + \phi}$)-regularity condition. To prove this starting from the definition, we consider the limits of tangent spaces 
to the fibres of $\psi + \phi$ along a sequence of points $(x_i, y_i) \in \bR^m \times \bR^n$
tending to some point $(\alpha, \beta) \in  W_1 \times W_2$ for some strata $W_1\in \cW_1$ and $W_2\in \cW_2$, at least one of which is of positive dimension.
We work with sequences such that $\psi(x_i) +\phi(y_i) \not=0$ for any $i >0$.

If $\psi(x_i) \not= 0$ and $\phi(y_i)\not= 0$ for any high enough index $i$, then we have the inclusion: $T_{(x_i, y_i)}(\psi + \phi) \supset T_{x_i} \psi \times T_{y_i} \phi$. Else, if say $\psi(x_i) =0$ for all $i$ high enough, then $\phi(y_i) \not=0$ for those indices $i$, and we have the inclusion: $T_{(x_i, y_i)}(\psi + \phi) \supset T_{x_i} W_{x_i} \times T_{y_i} \phi$, where $W_{x_i}$ denotes the stratum of $V(\psi)$ such that $x_i \in W_{x_i}$. 
Proving this inclusion essentially amounts to checking that for any $v\in T_{x_i} W_{x_i}$ one has $D \psi (x_i) v = 0$. Indeed,  if one considers a path $\alpha(t)$ within $W_{x_i}$ with $\frac{\partial \alpha}{\partial t}_{|t=0} = v$ then we have $\psi \circ \alpha \equiv 0$ thus $\frac{\partial (\psi \circ \alpha)}{\partial t}_{|t=0} =0$  and the latter is also equal to $D \psi (x_i) v$.

In both situations, taking limits we get the desired inclusion $\lim_{i\to \infty}T_{(x_i, y_i)}(\psi + \phi) \supset W_1 \times W_2$, due to the Thom regularity of $\psi$ and of $\phi$.

 Finally, one may refine the product stratification $\cW_1 \times \cW_2$ to a Whitney (a)-regular stratification such that the singular locus $\Sing(\psi + \phi)$ is a union of strata. This is possible by the classical theory of Whitney stratifications, see e.g. \cite{GWPL} for the algorithm. By construction, this refinement is 
 a Thom regular stratification of $V(\psi + \phi)$. This finishes the proof of our claim. 

 Taking now into account the codimension $p$ condition of the statement too (not used up to now), it appears that the map $\psi + \phi$ verifies the hypotheses of Theorem \ref{t:milnor} in view of our above remarks \S \ref{ss:thom} about the Thom condition. The  first claim of our statement follows.

Furthermore, if each mapping is radial weighted-homogeneous, the separate variable sum $\psi + \phi$ has the same property (of course, the weights have to be multiplied by some integers if the weighted degrees of $\psi$ and $\phi$ are not the same).   Thus one may apply Corollary \ref{c:radial}  and get the second claim of the statement.
\end{proof}

\begin{example}\label{e:ex2}
Let $h : \bR^3 \times \bC^n \to \bR^2$, $h= f(x,y,z) + g(w_1, \ldots , w_n)$, where  $f: \bR^3  \to \bR^2$ is Example \ref{e:thom} and $g : \bC^n \to \bC = \bR^2$ is a sum of monomials $\sum_{i=1}^{n} a_i m_i$, where $m_i = w_i^{k_i}$ or $m_i = \bar w_i^{k_i}$  with complex coefficients $a_i$. Then $g$ is radial weighted-homogeneous and has Thom property by Proposition \ref{p:exempl} or by the fact that this map has isolated singularity. As for $f$, it is not radial weighted-homogeneous, we have seen above that it satisfies the condition $\overline{M(f)\m V} \cap  V = \{ 0\}$ but we need the Thom regularity condition in order to apply Proposition \ref{p:exempl} (first claim). We claim that $f$ also satisfies the Thom $(a_{f})$-regularity condition at $\Sing f \m \{ 0\}$.
So let us fix some point $(0,0,z_{0})\in \Sing f= V(f)$, $z_{0} \not= 0$, and
choose an analytic curve $\gamma(t) =(x(t),y(t),z(t))$ with image in $B_{\varepsilon}\setminus V(f)$ 
defined on a small enough interval $\left]0,\varepsilon\right[$ such that $\lim_{t\to 0}\gamma(t) = (0,0,z_{0})$. For any $t$, the normal vector field $v_{1}(t)=(y(t),x(t),0)$ to the fibres of $f$ is orthogonal to the direction $(0,0,1)$ of the line $V(f)$. Let us consider the normal vector field to the fibres of $f$,  
 $v_{2}(t)=(-4x^{3}(t)-2x(t)z^{2}(t),4y^{3}(t), -2z(t)x^{2}(t))$, which is independent of $v_1(t)$. What we actually need to show is that the limit of the vector
product $\frac{v_{1}(t)}{\left\Vert v_{1}(t)\right\Vert }\wedge\frac{v_{2}(t)}{\left\Vert v_{2}(t)\right\Vert }$
is a nonzero vector in the space spanned by $(0,0,1)$.

One has $v(t):=v_{1}(t)\wedge v_{2}(t)=(-2x^{3}(t)z(t),2z(t)x^{2}(t)y(t),4y^{4}(t)+2x^{2}(t)z^{2}(t)+4x^{4}(t))$.  
 We have the following two situations:
\begin{itemize}
\item[(i)] If $\mathrm{ord}_{t}(x(t))\leq 2\mathrm{ord}_{t}(y(t))$  or if  $y(t) \equiv 0$,
we have:
\[
\lim_{t\rightarrow0}\frac{4y^{4}(t)+2x^{2}(t)z^{2}(t)+4x^{4}(t)}{x^{2}(t)} 
 =2z_{0}^{2}+\lim_{t\rightarrow0}\frac{4y^{4}(t)}{x^{2}(t)} > 0
\]
 and $\lim_{t\rightarrow0}(-2x(t)z(t),2z(t)y(t))=(0,0)$. 
\item[(ii)] If $\mathrm{ord}_{t}(x(t))>2 \mathrm{ord}_{t}(y(t))$ or if  $x(t) \equiv 0$, 
we have:\[
\lim_{t\rightarrow0}\frac{4y^{4}(t)+2x^{2}(t)z^{2}(t)+4x^{4}(t)}{y^{4}(t)}=4\neq0\]
 and $\lim_{t\rightarrow0}\frac{(-2x^{3}(t)z(t),2z(t)x^{2}(t)y(t))}{y^{2}(t)}=(0,0)$.
\end{itemize}

 These show that the tangent space to the line $\Sing f$ is contained in the limit of the tangent spaces to the fibres of $f$, hence that $f$ satisfies the Thom $(a_{f})$-regularity condition at $V(f) \m \{ 0\}$.

\end{example}
 
\subsection{Added in proof}
We have seen above that the Thom regularity condition along $V$ is not necessary for having the Milnor fibration in a conical neighbourhood of $V$, as explained in the first lines of \S \ref{s:strong}. Indeed, consider for instance the mixed function $f:\bC^2 \to \bC$, $f(x,y) = xy\bar x$. It is polar homogeneous, hence verifies condition \eqref{eq:main} and our Theorem \ref{t:homogen}, but one can easily check that there is no Thom (a$_f$) stratification of $V$. A deformation of this example in 3 variables has been suggested by A. Parusi\' nski during a recent workshop in Oberwolfach \cite{Top}, namely $g(x,y,z)= (x+z^3)y\bar x$. This is polar weighted-homogeneous and verifies condition \eqref{eq:main} and our Theorem \ref{t:homogen}, but again there is no Thom (a$_f$) stratification of $V$. This is a positive answer to our conjecture formulated above in the second paragraph of \S \ref{s:classicandnew}, i.e. that Thom's regularity condition is strictly stronger than condition \eqref{eq:main}. It 
is also a simple counter-example to a statement conjectured by Pichon and Seade in \cite[pag. 494]{PS} and proved as main result in a subsequent paper of the same authors.


\end{document}